\documentclass[12pt]{article}%
\usepackage{amsmath}
\usepackage{amsfonts}
\usepackage{amssymb}
\usepackage{graphicx}%
\makeatletter
\newfont{\footsc}{cmcsc10 at 8truept}
\newfont{\footbf}{cmbx10 at 8truept}
\newfont{\footrm}{cmr10 at 10truept}
\newtheorem{theorem}{Theorem}

\newtheorem{conjecture}[theorem]{Conjecture}

\newtheorem{proposition}[theorem]{Proposition}

\newenvironment{proof}[1][Proof]{\noindent{\textbf {#1}  }}  {\hfill$\blacksquare$\bigskip}

\def\blfootnote{\xdef\@thefnmark{}\@footnotetext}
\begin{document}

\title{\textbf{Maxima of the }$Q$\textbf{-index: graphs with no }$K_{s,t}$}
\author{Maria Aguieiras A. de Freitas\thanks{Federal University of Rio de Janeiro,
Brazil; \textit{email: maguieiras@im.ufrj.br} } , Vladimir
Nikiforov\thanks{Department of Mathematical Sciences, University of Memphis,
Memphis TN 38152, USA; \textit{email: vnikifrv@memphis.edu} } \ and Laura
Patuzzi\thanks{Federal University of Rio de Janeiro, Brazil; \textit{email:
laura@im.ufrj.br}} }
\date{}
\maketitle

\begin{abstract}
This note presents a new spectral version of the graph Zarankiewicz problem:
How large can be the maximum eigenvalue of the signless Laplacian of a graph
of order $n$ that does not contain a specified complete bipartite subgraph. A
conjecture is stated about general complete bipartite graphs, which is proved
for infinitely many cases.

More precisely, it is shown that if $G$ is a graph of order $n,$ with no
subgraph isomorphic to $K_{2,s+1},$ then the largest eigenvalue $q(G)$ of the
signless Laplacian of $G$ satisfies
\[
q(G)\leq\frac{n+2s}{2}+\frac{1}{2}\sqrt{(n-2s)^{2}+8s},
\]
with equality holding if and only if $G$ is a join of $K_{1}$ and an
$s$-regular graph of order $n-1.$\medskip



\end{abstract}

\section{Introduction.}

\emph{How many edges can a graph of order }$n$\emph{ have if it does not
contain a given complete bipartite subgraph? }This variant of the famous
Zarankiewicz problem \cite{Zar51} has turned out to be one of the most
difficult problems in modern Discrete Mathematics, widely open despite long
and intensive research. A comprehensive account of this vast theory can be
found in the survey of F\"{u}redi and Simonovits \cite{FuSi14}.

For our presentation, let us write $e\left(  G\right)  $ for the number of
edges of a graph $G$ and $K_{s,t}$ for the complete bipartite graph with
vertex classes of sizes $s$ and $t.$ Thus, the above problem can be stated
as:\medskip

\textbf{Problem A.} \emph{What is the maximum }$e\left(  G\right)  $\emph{ if
}$G$\emph{ is a graph of order }$n$\emph{ containing no }$K_{s,t}?\medskip$

Except for very few pairs of $s$ and $t$, no general solution of Problem A is
known. In a nutshell, the crucial difficulty is in the lack of constructions
proving that the known upper bounds on $e\left(  G\right)  $ are tight.

Improving the long standing results of K\"{o}vari, S\'{o}s, and Tur\'{a}n
\cite{KST54} and of Zn\'{a}m \cite{Zna63}, F\"{u}redi \cite{Fur96} gave a
general bound, which was polished by Nikiforov in \cite{Nik10b} to the final
form stated in Theorem \ref{tFu} below.

In fact, Nikiforov considered a spectral version of Problem A, which was also
studied by Babai and Guiduli in \cite{BaGu09}. Thus, let $\lambda\left(
G\right)  $ denote the spectral radius of the adjacency matrix of a graph $G.$
The following problem is the natural spectral analog of Problem A:\medskip\ 

\textbf{Problem B. }\emph{What is the maximum }$\lambda\left(  G\right)  $
\emph{of a graph }$G$ \emph{of order }$n$ \emph{containing no }$K_{s,t}%
?\medskip$

An asymptotic upper bound was given by Babai and Guiduli in \cite{BaGu09};
independently the following more precise statement was proved in \cite{Nik07d}
and \cite{Nik10b}:

\begin{theorem}
\label{th1}Let $s\geq t\geq2,$ and let $G$ be a $K_{s,t}$-free graph of order
$n,$ with spectral radius $\lambda.$ If $t=2,$ then%
\begin{equation}
\lambda\leq1/2+\sqrt{\left(  s-1\right)  \left(  n-1\right)  +1/4}.
\label{in0}%
\end{equation}
If $t\geq3,$ then%
\begin{equation}
\lambda\leq\left(  s-t+1\right)  ^{1/t}n^{1-1/t}+\left(  t-1\right)
n^{1-2/t}+t-2. \label{in1}%
\end{equation}

\end{theorem}

Note that Theorem \ref{th1} is closely related to Problem A: indeed, the
well-known inequality $\lambda\left(  G\right)  \geq2e\left(  G\right)  /n$
immediately implies F\"{u}redi's result \cite{Fur96}, and even a slight
improvement of it.

\begin{theorem}
\label{tFu}Let $s\geq t\geq2,$ and let $G$ be a $K_{s,t}$-free graph of order
$n,$ with $e\left(  G\right)  $ edges$.$ If $t=2,$ then%
\begin{equation}
e\left(  G\right)  \leq\frac{n}{2}\sqrt{\left(  s-1\right)  \left(
n-1\right)  +1/4}+\frac{n}{4}. \label{HC}%
\end{equation}
If $t\geq3,$ then%
\[
e\left(  G\right)  \leq\frac{1}{2}\left(  s-t+1\right)  ^{1/t}n^{2-1/t}%
+\frac{1}{2}\left(  t-1\right)  n^{2-2/t}+\frac{1}{2}\left(  t-2\right)  n.
\]

\end{theorem}

It is worth pointing out that inequality (\ref{HC}) follows from a theorem of
Hylt\'{e}n-Cavallius \cite{HyCa58} and is one of the few known tight results
in the area, since F\"{u}redi constructed a matching family of graphs in
\cite{Fur96a}.

\begin{theorem}
For any $n$ there exist a $K_{2,s+1}$-free graph $G$ of order $n$ such that%
\[
e\left(  G\right)  \geq\frac{1}{2}n\sqrt{sn}+O\left(  n^{4/3}\right)  ,
\]

\end{theorem}

Another point to be made here is that equality holds in (\ref{in0}) and
(\ref{HC}) if and only if $G$ is a strongly regular graph with parameters
\[
\left(  n,1/2+\sqrt{\left(  s-1\right)  \left(  n-1\right)  +1/4},s,s\right)
.
\]
Such strongly regular graphs are sometimes called \emph{design graphs }(see,
e.g., \cite{BJL00}) and appear in various problems.\medskip

The similarity between Theorems \ref{th1} and \ref{tFu}, together with the
fact that most known constructions for Problem A are regular or almost regular
graphs, suggests that Problems A and B might be essentially equivalent, and
therefore equally hard. Let us note that such equivalence has been indeed
proved for nonbipartite forbidden subgraphs. Thus, it is of interest whether
other spectral versions of Problem A may be substantially different from
Problem B, like the following one:\medskip

\textbf{Problem C. }\emph{What is the maximum spectral radius of the signless
Laplacian of a graph }$G$ \emph{of order }$n$ \emph{containing no }%
$K_{s,t}?\medskip$

To make this question clear we need a brief introduction: let $G$ be a graph
with adjacency matrix $A$ and let $D$ be the diagonal matrix of the row-sums
of $A$, i.e., the degrees of $G$. The matrix $Q(G)=A+D$, called the
\emph{signless Laplacian} or the $Q$-matrix of $G$, has been intensively
studied, see, e.g., the survey of Cvetkovi\'{c} \cite{C10} and its references.
The maximal eigenvalue (equivalently, the spectral radius) of $Q(G)$ is called
the $Q$\emph{-index} of $G$ and is denoted by $q\left(  G\right)  $.

In general, in extremal problems with forbidden nonbipartite graphs,
$\lambda\left(  G\right)  $ and $q\left(  G\right)  $ behave similarly (see a
discussion of this fact in \cite{FNP13}), but they may be considerable
differences between extremal problems about $\lambda\left(  G\right)  $ and
$q\left(  G\right)  $ in case of forbidden bipartite graphs.

Regarding Problem C the authors believe it is not as difficult as Problems A
and B, and it will be resolved completely in the next few years. Moreover, we
venture the following conjecture:

\begin{conjecture}
\label{con}Let $s\geq t-1\geq1$, and let $n$ be sufficiently large. If $G$ is
a $K_{t,s+1}$-free graph of order $n$, then
\[
q\left(  G\right)  \leq\frac{n}{2}+s+t-2+\frac{1}{2}\sqrt{\left(
n-2+2s\right)  ^{2}-8s\left(  n-2\right)  +4\left(  t-1\right)  \left(
n-t+1\right)  }.
\]
Equality holds if and only if $G$ is a join of $K_{t-1}$ and an $s$-regular
graph of order $n-t+1$.
\end{conjecture}

Presently we cannot prove or disprove this conjecture for all $s$ and $t.$ For
$t=2$ and $s=1$ the graph $K_{2,2}$ is just the cycle of length $4,$ and this
case of Conjecture \ref{con} was confirmed in \cite{FNP13}. In this paper we
shall solve the case $t=2$ and any $s\geq1.$

\begin{theorem}
\label{thc1}Let $s\geq1$, and let $n\geq s^{2}+6s+6.$ If $G$ is a $K_{2,s+1}$
graph of order $n,$ containing no $K_{2,s+1}$, then
\begin{equation}
q\left(  G\right)  \leq\dfrac{n+2s}{2}+\dfrac{1}{2}\sqrt{(n-2s)^{2}+8s}.
\label{mmi}%
\end{equation}
Equality holds if and only if $G$ is a join of $K_{1}$ and an $s$-regular
graph of order $n-1$.
\end{theorem}

We shall break Theorem \ref{thc1} into two separate statements, with separate
proofs. The purpose of this separation is twofold: first it helps with the
presentation of the proof, and second it may be easier to analyze the proof
and extend it for the general case of Conjecture \ref{con}.\medskip

We start with a result about a join of a vertex with a graph of bounded
maximum degree. As proved in \cite{MNRS}, if $G$ is a join of $K_{1}$ and an
$s$-regular graph of order $n-1,$ then%
\[
q\left(  G\right)  =\dfrac{n+2s}{2}+\dfrac{1}{2}\sqrt{(n-2s)^{2}+8s}.
\]
We shall give an easy improvement of this statement, making it an extremal
result. Write $G\vee H$ for the join of two graphs $H$ and $G.$

\begin{proposition}
\label{pro}Let $s\geq1$, let $H$ be a graph of order $n-1,$ and let
$G:=K_{1}\vee H.$ If $\Delta\left(  H\right)  \leq s,$ then
\begin{equation}
q\left(  G\right)  \leq\dfrac{n+2s}{2}+\dfrac{1}{2}\sqrt{(n-2s)^{2}+8s}
\label{mi}%
\end{equation}
Equality holds if and only if $H$ is $s$-regular.
\end{proposition}

We postpone the proof of Proposition \ref{pro} to Section \ref{pfs}, after we
introduce the necessary notation. Here we just note that simple as it is,
Proposition \ref{pro} immediately takes care of the essential case of Theorem
\ref{thc1}.

\begin{theorem}
\label{th2}Let $s\geq1$, and let $G$ be a graph of order $n,$ with
$\Delta\left(  G\right)  =n-1.$ If $G$ is $K_{2,s+1}$-free, then
\[
q\left(  G\right)  \leq\dfrac{n+2s}{2}+\dfrac{1}{2}\sqrt{(n-2s)^{2}+8s}.
\]
Equality holds if and only if $H$ is an $s$-regular graph.
\end{theorem}

Theorem \ref{th2} obviously follows from Proposition \ref{pro}, so we shall
omit its proof.

The following theorem completes the proof of Theorem \ref{thc1} and shows that
if the premise $\Delta\left(  G\right)  =n-1$ is relaxed, we can strengthen
the bound on $q\left(  G\right)  $.

\begin{theorem}
\label{th3}Let $s\geq1$, $n\geq s^{2}+6s+6,$ and $G$ be a graph of order $n,$
with $\Delta\left(  G\right)  <n-1.$ If $G$ is $K_{2,s+1}$-free, then
$q\left(  G\right)  <n.$
\end{theorem}

Much of the rest of the paper is dedicated to the proof of Theorem \ref{th3},
which is not too short.

\section{\label{pfs}Proofs of Proposition \ref{pro} and Theorem \ref{th3}}

First we shall introduce some notation; for graph notation undefined here we
refer the reader to \cite{Bol98}. Thus, if $G$ is a graph, and $X$ and $Y$ are
disjoint sets of vertices of $G$, we write:\medskip

- $V\left(  G\right)  $ for the set of vertices of $G$;

- $E(G)$ for the set of edges of $G,$ and let $e\left(  G\right)  :=|E(G)|$;

- $\Delta\left(  G\right)  $ for the maximum degree of $G$;

- $\Gamma\left(  u\right)  $ for the set of neighbors of a vertex $u$, and let
$d(u):=\left\vert \Gamma\left(  u\right)  \right\vert $;

- $G\left[  X\right]  $ for the graph induced by $X$, and let $E\left(
X\right)  :=E\left(  G\left[  X\right]  \right)  $ and $e\left(  X\right)
:=|E(G)|$;

- $e\left(  X,Y\right)  $ for the number of edges joining vertices in $X$ to
vertices in $Y.\medskip$

Regarding the right side of (\ref{mmi}), note that
\[
\dfrac{n+2s}{2}+\dfrac{1}{2}\sqrt{(n-2s)^{2}+8s}>n,
\]
and also
\begin{align}
\dfrac{n+2s}{2}+\dfrac{1}{2}\sqrt{(n-2s)^{2}+8s}  &  =n+\dfrac{\sqrt
{(n-2s)^{2}+8s}-\left(  n-2s\right)  }{2}\nonumber\\
&  =n+\dfrac{4s}{\sqrt{(n-2s)^{2}+8s}+\left(  n-2s\right)  }<n+\dfrac
{2s}{n-2s}. \label{bo}%
\end{align}
\medskip

Part of our proof of Theorem \ref{th3} is based of the following inequality
that can be traced back to Merris \cite{Mer98}:\medskip

\emph{For every graph }$G,$
\begin{equation}
q\left(  G\right)  \leq\max\left\{  d\left(  u\right)  +\frac{1}{d\left(
u\right)  }\sum_{v\in\Gamma\left(  u\right)  }d\left(  v\right)  \text{
$\vert$
}u\in V\left(  G\right)  \text{ and }d\left(  u\right)  >0\right\}  .
\label{Mb}%
\end{equation}
\medskip

\begin{proof}
[\textbf{Proof of Proposition \ref{pro}}]Let $w$ be the vertex of $G$
corresponding to $K_{1}$ in the representation $G:=K_{1}\vee H.$ Set for short
$q:=q\left(  G\right)  $ and let $\mathbf{x}:=\left(  x_{1},\ldots
,x_{n}\right)  $ be a positive eigenvector to $q.$ Choose a vertex $u\in
V\left(  H\right)  $ such that
\[
x_{u}=\max_{a\in V\left(  H\right)  }x_{a}%
\]
From the eigenequations for the $Q$-matrix we have
\begin{equation}
qx_{w}=\left(  n-1\right)  x_{w}+\sum_{v\in V\left(  H\right)  }x_{v}%
\leq\left(  n-1\right)  x_{w}+(n-1)x_{u} \label{i0}%
\end{equation}
and%
\[
qx_{u}=d(u)x_{u}+\sum_{\left\{  v,u\right\}  \in E\left(  G\right)  }x_{v}%
\leq\left(  s+1\right)  x_{u}+x_{w}+sx_{u}.
\]
Hence, we find that%
\begin{align}
\left(  q-n+1\right)  x_{w}  &  \leq\left(  n-1\right)  x_{u},\label{i1}\\
\left(  q-2s-1\right)  x_{u}  &  \leq x_{w}. \label{i2}%
\end{align}
On the other hand, it is known that $q\geq\Delta\left(  G\right)  +1$; thus,
$q-n+1>0.$ Therefore, we can multiply (\ref{i1}) and (\ref{i2}), obtaining%
\[
q^{2}-\left(  n+2s\right)  q-2s\left(  n-1\right)  \leq0,
\]
which implies (\ref{mi}).

If equality holds in (\ref{mi}), then equality holds in (\ref{i0}), and so
$x_{v}=x_{u}$ for any vertex $v\in V\left(  H\right)  .$ Since for any $v\in
V\left(  H\right)  $ we have%
\[
qx_{v}=d(v)x_{v}+\sum_{\left\{  p,v\right\}  \in E\left(  G\right)  }x_{p}%
\leq\left(  s+1\right)  x_{u}+x_{w}+sx_{u}=qx_{u}%
\]
we see that $d(v)=s+1,$ and so $H$ is $s$-regular.
\end{proof}

\medskip

\begin{proof}
[\textbf{Proof of Theorem \ref{th3}}]Suppose that $G$ satisfies the hypothesis
of the theorem, and for any nonisolated vertex $u\in V\left(  G\right)  ,$
let
\[
F\left(  u\right)  :=d\left(  u\right)  +\frac{1}{d\left(  u\right)  }%
\sum_{v\in\Gamma\left(  u\right)  }d\left(  v\right)  .
\]
Our first goal is to prove that if $\Delta\left(  G\right)  \leq n-s-2,$ then
$F\left(  u\right)  <n$ for any nonisolated $u\in V\left(  G\right)  ,$ which,
in view of (\ref{Mb}), implies that $q\left(  G\right)  <n$ as well$.$

If $d(u)\leq s+1$, then $F\left(  u\right)  <n$ follows immediately by
\[
d(u)+\frac{1}{d(u)}\sum_{s\in\Gamma(u)}d(s)\leq s+1+\Delta(G)<n,
\]
so from now on we shall suppose that $s+2\leq d(u)\leq n-s-2.$

Fix a nonisolated vertex $u\in V\left(  G\right)  ,$ and let $A:=\Gamma\left(
u\right)  $ and $B:=V\left(  G\right)  \backslash\left(  A\cup\left\{
u\right\}  \right)  .$ We see that
\[
F\left(  u\right)  =d(u)+1+\frac{1}{d(u)}\left(  2e\left(  A\right)  +e\left(
A,B\right)  \right)  .
\]
Since $G$ is $K_{2,s+1}$-free, we have $\Delta\left(  G\left[  A\right]
\right)  \leq s$ and $\left\vert \Gamma\left(  v\right)  \cap A\right\vert
\leq s$ for any $v\in B.$ Hence,
\[
2e\left(  A\right)  \leq d(u)s
\]
and
\[
e\left(  A,B\right)  =\sum_{v\in B}\left\vert \Gamma\left(  v\right)  \cap
A\right\vert \leq\sum_{v\in B}s=(n-d(u)-1)s.
\]
Adding the last two inequalities, we get the bound
\[
F\left(  u\right)  \leq d(u)+1+\frac{(n-1)s}{d(u)},
\]
so we want to prove that the right side is always less than $n.$ Since the
function
\[
g(x):=x+1+\frac{(n-1)s}{x}%
\]
is convex for $x>0$, and $s+2\leq d(u)\leq n-s-2$, we see that
\begin{align*}
F\left(  u\right)   &  \leq\max\{g(s+2),g(n-s-2)\}\\
&  =\max\left\{  s+3+\frac{(n-1)s}{s+2},n-s-1+\frac{(n-1)s}{n-s-2}\right\}
\end{align*}
After some simple algebra, we find that if $n>s^{2}+2s+2,$ then%
\[
\max\left\{  s+3+\frac{(n-1)s}{s+2},n-s-1+\frac{(n-1)s}{n-s-2}\right\}  <n,
\]
completing the proof if $\Delta\left(  G\right)  \leq n-s-2.$

Thus, it remains to prove the theorem for $\Delta\left(  G\right)  \geq
n-s-1.$ In this case we shall use completely different approach. Assume for a
contradiction that $q\left(  G\right)  >n.$ Set for short $q:=q\left(
G\right)  $ and let $\mathbf{x}:=\left(  x_{1},\ldots,x_{n}\right)  $ be a
nonnegative unit eigenvector to $q$.

Choose a vertex $w$ with $d(w)=\Delta(G)$ and let $A=\Gamma(w)$ and
$B=V(G)\setminus\left(  \{w\}\cup A)\right)  $. Note that
\[
\left\vert B\right\vert =n-1-d\left(  w\right)  =n-1-\Delta(G)\leq s.
\]
\qquad Now, choose $u\in A$ and $v\in B$ such that%
\[
x_{u}=\max_{a\in A}x_{a}\text{ \ \ and \ \ \ }x_{v}=\max_{a\in B}x_{a}.
\]
Our first goal is to show that $d(u)\leq2s+1$ and $d(v)\leq2s-1.$ This follows
easily by
\[
d(u)=1+\left\vert \Gamma\left(  u\right)  \cap A\right\vert +\left\vert
\Gamma\left(  u\right)  \cap B\right\vert \leq1+s+\left\vert B\right\vert
\leq2s+1,
\]
and likewise,%
\[
d(v)=\left\vert \Gamma(v)\cap B\right\vert +\left\vert \Gamma(v)\cap
A\right\vert \leq\left\vert B\right\vert -1+s\leq2s-1.
\]

Next, we shall show that $x_{v}<2/n$ and $x_{u}<2/n.$ To this end note that
the eigenequations for $Q\left(  G\right)  $ corresponding to $u$ and $v$
imply that%
\[
qx_{u}=d\left(  u\right)  x_{u}+x_{w}+\sum_{a\in\Gamma(u)\cap A}x_{a}%
+\sum_{a\in\Gamma(u)\cap B}x_{a}<\left(  2s+1\right)  x_{u}+1+sx_{u}+sx_{v},
\]
and
\[
qx_{v}=d\left(  v\right)  x_{v}+\sum_{a\in\Gamma(v)\cap A}x_{a}+\sum
_{a\in\Gamma(v)\cap B}x_{a}\leq\left(  2s-1\right)  x_{v}+sx_{u}+\left(
s-1\right)  x_{v}.
\]
Rearranging these inequalities, we get%
\begin{equation}
\left(  q-3s-1\right)  x_{u}<1+sx_{v} \label{i3}%
\end{equation}
and%
\[
\left(  q-3s+2\right)  x_{v}\leq sx_{u}.
\]
Now, excluding $x_{u}$, we find that
\[
\left(  q-3s+2\right)  x_{v}<\frac{\left(  1+sx_{v}\right)  s}{q-3s-1}%
\leq\frac{\left(  1+sx_{v}\right)  s}{n-3s-1}<\frac{s}{n-3s-1}+sx_{v},
\]
and so,%
\[
x_{v}<\frac{s}{\left(  n-3s-1\right)  \left(  n-4s+2\right)  }\leq\frac
{1}{\left(  s+1\right)  \left(  n-4s+2\right)  }\leq\frac{1}{2\left(
n-4s+2\right)  }\leq\frac{2}{n}.
\]
To bound $x_{u}$ we substitute $2/n$ for $x_{v}$ in (\ref{i3}) and obtain%
\[
x_{u}<\frac{1+sx_{v}}{q-3s-1}<\frac{1+2s/n}{n-3s-1}<\frac{1+1/3}{n-3s-1}%
\leq\frac{2}{n}.
\]

Armed with the upper bounds on $x_{u}$ and $x_{v},$ we shall prove that $q<n.$
We shall use some relatively new techniques for this purpose. To begin with,
since $\mathbf{x}$ is a unit vector, we have%
\[
q=\left\langle Q\mathbf{x},\mathbf{x}\right\rangle =\sum\limits_{\left\{
i,j\right\}  \in E(G)}\left(  x_{i}+x_{j}\right)  ^{2}.
\]
Write $G^{\prime}$ for the graph $G\left[  A\cup\left\{  w\right\}  \right]  $
and set $n^{\prime}=\left\vert V\left(  G^{\prime}\right)  \right\vert .$ Note
that%
\begin{equation}
\sum\limits_{\left\{  i,j\right\}  \in E(G)}\left(  x_{i}+x_{j}\right)
^{2}=\sum\limits_{\left\{  i,j\right\}  \in E(G^{\prime})}\left(  x_{i}%
+x_{j}\right)  ^{2}+\sum\limits_{\left\{  i,j\right\}  \in E(B)}\left(
x_{i}+x_{j}\right)  ^{2}+\sum\limits_{\left\{  i,j\right\}  \in E(A,B)}\left(
x_{i}+x_{j}\right)  ^{2} \label{dec}%
\end{equation}
Since $\Delta\left(  G^{\prime}\right)  =n^{\prime}-1,$ Theorem \ref{th2},
together with (\ref{bo}), implies that
\[
\sum\limits_{\left\{  r,s\right\}  \in E(G^{\prime})}(x_{r}+x_{s})^{2}%
\leq\dfrac{n^{\prime}+2s+\sqrt{(n^{\prime}-2s)^{2}+8s}}{2}<n^{\prime}%
+\dfrac{2s}{n^{\prime}-2s}.
\]
On the other hand, using the inequalities $x_{v}<2/n$ and $x_{u}<2/n,$ we see
that
\begin{align*}
\sum\limits_{\left\{  i,j\right\}  \in E(B)}\left(  x_{i}+x_{j}\right)
^{2}+\sum\limits_{\left\{  i,j\right\}  \in E(A,B)}\left(  x_{i}+x_{j}\right)
^{2}  &  \leq\sum\limits_{\left\{  i,j\right\}  \in E(B)}(x_{v}+x_{v}%
)^{2}+\sum\limits_{\left\{  i,j\right\}  \in E(A,B)}(x_{u}+x_{v})^{2}\\
&  \leq e\left(  B\right)  \left(  \frac{2}{n}+\frac{2}{n}\right)
^{2}+e\left(  A,B\right)  \left(  \frac{2}{n}+\frac{2}{n}\right)  ^{2}\\
&  <\left(  \binom{s}{2}+s^{2}\right)  \frac{16}{n^{2}}<\frac{24s^{2}}{n^{2}}.
\end{align*}
Therefore, in view of (\ref{bo}) and (\ref{dec}), we obtain
\[
q\leq n^{\prime}+\dfrac{2s}{n^{\prime}-2s}+\frac{24s^{2}}{n^{2}}.
\]
Since the function $g\left(  x\right)  =x+2s/\left(  x-2s\right)  $ is convex
whenever $x>2s,$ the inequalities
\[
n-s\leq n^{\prime}\leq n-1
\]
imply that
\[
n^{\prime}+\dfrac{2s}{n^{\prime}-2s}\leq\max\left\{  n-1+\frac{2s}%
{n-1-2s},n-s+\frac{2s}{n-3s}\right\}  .
\]
In view of $s\geq1$ and $n>3s+2,$ one easily finds that
\[
n-s+\frac{2s}{n-3s}\leq n-1+\frac{2s}{n-1-2s}.
\]
Therefore,%
\begin{align*}
q  &  \leq n-1+\frac{2s}{n-1-2s}+\frac{24s^{2}}{n^{2}}\leq n-1+\frac
{2}{s+4+5/s}+\frac{24}{\left(  s+6+6/s\right)  ^{2}}\\
&  <n-1+\frac{2}{5}+\frac{24}{49}<n.
\end{align*}
The proof of Theorem \ref{th3} is completed.
\end{proof}

\section{Concluding remarks}

In our proof of Theorems \ref{thc1} and \ref{th3} we used techniques that have
worked efficiently for solving a number of extremal problems about the
$Q$-index; however, these methods seem inadequate for tackling Conjecture
\ref{con} in general. We need completely new general techniques, for which
Conjecture \ref{con} provides both motivation and a test field.\medskip

\textbf{Acknowledgement}\medskip

The first author was partially supported by CNPq (the Brazilian Council for
Scientific and Technological Development) and FAPERJ (Foundation for Research
of the State of Rio de Janeiro).

\end{document}